# Trigonometric splines of several variables

Denysiuk V.P.

Dr. of Phys-Math. sciences, Professor, Kiev, Ukraine

National Aviation University

kvomden@nau.edu.ua

## Annotation

Under consideration methods of constructing trigonometric interpolation splines of two variables on rectangular areas. These methods are easily generalized to the case of trigonometric interpolation splines of several variables on such domains. A numerical example illustrating the main theoretical propositions is considered. The given methods of constructing such splines can be widely used in practice.

## Keywords:

## Introduction

Polynomial splines of two variables are piecewise polynomials that are stitched together according to a certain pattern. They are of great interest in the theory of approximations, numerical analysis, construction and reconstruction of various surfaces in applied fields, when solving boundary problems of mathematical physics by numerical methods, etc.

Probably, the first systematic presentation of the theory of polynomial splines in one or two variables was carried out in [1]. Methods of constructing such splines were also considered there.

Polynomial splines of two variables have been intensively studied in recent years. The main results obtained from these studies up to the year 2000 are discussed in [2]. Methods of interpolation by two-dimensional splines (including numerical examples) and generalization of results about dimensionality and the order of approximation by such splines are given there. It is also emphasized that, in contrast to the one-dimensional case, the solution of even standard problems for the two-dimensional case poses significant difficulties. In particular, the construction of explicit interpolation schemes (especially Lagrangian interpolation schemes) for two-dimensional splines on given triangulations leads to complex problems.

A voluminous bibliography of works devoted to the construction and application of polynomial two-dimensional splines is given in [3].

Recently, the author proposed a new class of functions - trigonometric splines [4-5], which are special rapidly converging trigonometric Fourier series; at the same time, it was possible to combine the developed theory of trigonometric series with the theory of polynomial splines.

The construction of trigonometric splines of several variables looks attractive. In this work, two methods of constructing trigonometric splines of two variables on rectangular regions are considered.

## The purpose of the work.

Development of methods for constructing trigonometric splines of two variables on rectangular areas.

# The main part.

Leave it for a while $[0.2\pi)$ a uniform grid is given $\Delta_N^{(I)} = \left\{ x_j^{(I)} \right\}_{j=1}^{N}$, where $I$ -grid indicator,

( $I = 0,1$ ), where $x_j^{(0)} = \dfrac{2\pi}{N}(j-1)$, $x_j^{(1)} = \dfrac{\pi}{N}(2j-1)$, $N = 2n+1$, $n = 1,2,\ ....$ Let it too $C_p^r$ ( $r = 0,1,2,\ ...$ ) is

space $2\pi$ -periodic, continuous functions having continuous derivatives up to order $r$ including; note that with this symbol $C_p^{-1}$ we denote the set $2\pi$ -periodic, piecewise constant functions with a finite number of discontinuities of the first kind. Consider the function $f(x) \in C^r$, which takes in grid nodes $\Delta_N^{(I)}$ value $\left\{ f(x_j^{(I)}) \right\}_{j=1}^{N} = \left\{ f_j^{(I)} \right\}_{j=1}^{N}$.

Interpolation spline $S(f, r, \Delta_N, x)$ we will call the function $S(f, r, \Delta_N, x)$, which satisfies the following conditions:

and) $S(f, \Delta_N, x_j^{(I)}) = f_j^{(I)}$, ( $j = 1,2,...,N$; $I = 0,1$ ) ;

in) $S(f, r, \Delta_N, x) \in C^r$, ( $r = 0,1,2,\ ...$ ).

Note that we have excluded from the traditional definition of splines the requirement that this spline coincide on each segment $\left[ x_{j-1}^{(0)}, x_j^{(0)} \right]$, ( $j = 2,3,...,N$ ), grids $\Delta_N^{(0)}$ with an algebraic polynomial of a certain degree, since this requirement actually determines the method of spline construction.

Currently, two methods of constructing interpolation splines are known $S(f, r, \Delta_N, x)$.

With one of them, this spline is constructed in the form of a set of algebraic polynomials that are given on each segment $\left[ x_{j-1}^{(0)}, x_j^{(0)} \right]$, ( $j = 2,3,...,N$ ), grids $\Delta_N^{(0)}$, and which are stitched together in such a way as to ensure the continuity of derivatives to $r$ of order inclusive. This method of constructing splines, which are called polynomial, was probably first studied in detail in[1]and became widespread. When constructing polynomial splines by this method, it is necessary to solve a system of algebraic equations with diagonal matrices. The theory of polynomial interpolation splines is well developed; so, in particular, it is known that on certain classes of smooth functions such splines provide the smallest approximation errors in some spaces[7].

At the same time, polynomial interpolation splines have a number of disadvantages, which, in our opinion, include the need to consider splines of even and odd degrees separately, the complexity of constructing splines of high degrees, and their piecemeal structure. These shortcomings lead to the fact that splines of the odd (most often third) degree are mostly used in practice; difficulties also arise when applying linear operators to polynomial splines.

The second construction method of interpolation splines $S(f, r, \Delta_N, t)$ consists in the fact that the spline is constructed in the form of a special rapidly converging trigonometric Fourier series; the convergence of this series is ensured by the convergence factors, the order of which determines the differential properties of the spline. This method of constructing splines, which can naturally be called trigonometric interpolation splines, was first given in [4,5]; further studies of this method are given in [6,8].

Trigonometric interpolation splines have a number of advantages compared to polynomial interpolation splines. So, for example, there is no need to distinguish between the methods of constructing trigonometric splines of even and odd degrees, since the construction of splines of arbitrary integer degrees is carried out in the same way. The application of linear operators to trigonometric splines also does not cause difficulties. It should be noted that at certain values of the parameters defining trigonometric splines, they coincide with periodic simple polynomial splines of both even and odd degrees; this allows us to extend to such trigonometric splines all the results of the theory of approximations obtained for polynomial splines. Instead of a system of algebraic equations that must be solved when constructing polynomial splines, when constructing trigonometric splines, it is necessary to find discrete coefficients of interpolating trigonometric polynomials; when finding these coefficients, it is possible (and expedient) to use the FFT (fast Fourier trans-form) method. It is easy to construct systems of fundamental trigonometric splines [9], which is important in many cases. Finally, the set of trigonometric interpolating splines is wider than the set of polynomial splines, since most trigonometric splines do not have polynomial counterparts.

Certain disadvantages of trigonometric splines include their periodicity; this requires the use of certain methods of constructing these splines when approximating non-periodic functions [10]. The presentation of trigonometric splines in the form of infinite trigonometric series can also be considered a disadvantage; however, since these series are rapidly converging, this drawback can be neglected.

Trigonometric interpolation splines of two variables can be presented in two ways. In one of them, a trigonometric interpolation polynomial of two variables is used, and trigonometric splines of two variables are constructed on the basis of this polynomial. The construction of trigonometric interpolation splines of several variables is reduced to the construction of the trigonometric interpolation polynomial of these variables.

In the second method, trigonometric splines of two variables are obtained through fundamental splines of two variables; in turn, fundamental splines of two variables are obtained as products of fundamental splines for each variable. This method of constructing trigonometric splines is simpler and easily generalizes to the case of trigonometric splines of several variables, which is useful in many applications.

Let's consider trigonometric splines of one variable in more detail.

Let us consider a trigonometric polynomial

$$T_n^{(I)}(x) = \frac{a_0^{*(I)}}{2} + \sum_{k=1}^{\frac{N-1}{2}} a_k^{(I)} \cos kx + b_k^{(I)} \sin kx, \quad (1)$$

which interpolates the function $f(t)$ on the grid $\Delta_N^{(I)}$; the coefficients of this polynomial are determined by formulas

$$a_k^{(I)} = \frac{2}{N} \sum_{j=1}^{N} f_j^{(I)} \cos kx_j^{(I)}, \quad b_k^{(I)} = \frac{2}{N} \sum_{j=1}^{N} f_j^{(I)} \sin kx_j^{(I)}, \quad (2)$$

$$k = 0,1,...,n; \qquad\qquad k = 1,2,...,n.$$

We will also introduce the concepts of stitching mesh and interpolation mesh. We will call the stitching mesh the mesh on which the polynomial analogy of the trigonometric splines, of course, in the case when they exist. An interpolation grid is a grid whose nodes are interpolated by the corresponding trigonometric polynomial $T_n^{(I)}(t)$. The stitching grid and the interpolation grid may or may not coincide with each other.

In [11] it was shown that the trigonometric spline with parameter vectors $\Gamma = \{\gamma_1, \gamma_2, \gamma_3\}$ and $H = \{\eta_1, \eta_2, \eta_3\}$, the components of which are arbitrary real numbers, and at least one of the components $\gamma_2, \gamma_3, \eta_2, \eta_3$ is not equal to 0, has the form .

$$St(I_1, I_2, \Gamma, H, \nu, r, t) = \frac{a_0^{(I_2)}}{2} + \sum_{k=1}^{\frac{N-1}{2}} \left[ a_\kappa^{(I_2)} Cs(I_1, I_2, \Gamma \nu, r, t) + b_k^{(I_2)} Ss(I_1, I_2, H, \nu, r, t) \right], \quad (3)$$

where

$$Cs_k(I_1, I_2, \Gamma \nu, r, t) = \frac{c_k(I_1, \Gamma, \nu, r, t)}{hc_k(I_1, I_2, \Gamma, \nu, r)} \qquad Ss_k(I_1, I_2, \Gamma \nu, r, t) = \frac{s_k(I_1, H, \nu, r, t)}{hs_k(I_1, I_2, H, \nu, r)} \quad (4)$$

$$c_k(I_1, \Gamma, \nu, r, t) = \gamma_1 \nu_k(r) \cos kt +$$
$$+ \sum_{m=1}^{\infty} (-1)^{ml_1} \left[ \gamma_3 \nu_{mN+k}(r) \cos((mN+k)t) + \gamma_2 \nu_{mN-k}(r) \cos((mN-k)t) \right]; \quad (5)$$

$$s_k(I_1, H, \nu, r, t) = \eta_1 \nu_k(r) \sin kt +$$
$$+ \sum_{m=1}^{\infty} (-1)^{ml_1} \left[ \eta_3 \nu_{mN+k}(r) \sin((mN+k)t) - \eta_2 \nu_{mN-k}(r) \sin((mN-k)t) \right];$$

with convergence factors $\nu_k(r) = (k)^{-(1+r)}$, and interpolation factors

$$hc_k(I_1, I_2, \Gamma, \nu, r) = \gamma_1 \nu_k(r) + \sum_{m=1}^{\infty} (-1)^{m(I_1-I_2)} \left[ \gamma_3 \nu_{mN+k}(r) + \gamma_2 \nu_{mN-k}(r) \right];$$

$$hs_k(I_1, I_2, H, \nu, r) = \eta_1 \nu_k(r) + \sum_{m=1}^{\infty} (-1)^{m(I_1-I_2)} \left[ \eta_3 \nu_{mN+k}(r) + \eta_2 \nu_{mN-k}(r) \right]. \quad (6)$$

where is the indicator $I_1$ $(I_1 = 0,1)$ defines the stitching grid, the indicator $I_2$ $(I_2 = 0,1)$ defines the interpolation grid, $a_0^{(I_2)}, a_\kappa^{(I_2)}, b_\kappa^{(I_2)}$ - coefficients of the interpolation trigonometric polynomial on the grid $\Delta_N^{(I_2)}$, $r$, $(r = 1, 2, ...)$ - a parameter that determines the order of the spline. Note that to reduce entries in the notation of splines and functions through which they are constructed, we omit the dependence on the number $N$ nodes of interpolation grids $\Delta_N^{(I)}$. In the following functions $Cs_k(I_1, I_2, \Gamma \nu, r, t)$ and $Ss_k(I_1, I_2, \Gamma \nu, r, t)$ will be called spline cos- and sin-functions, respectively.

It is clear that in the accepted notations $St(I_1, I_2, \Gamma, H, \nu, r, t) \in C_p^{r-1}$.

In the future, without losing generality, we will limit ourselves to considering only simple splines, when the parameter vectors have the form $\Gamma = \{1,1,1\}$ and $H = \{1,1,1\}$; in the notation of simple trigonometric splines, we will omit the dependence on vectors $\Gamma$ and $H$.

It is easy to see that expressions (1) and (3) are similar; An important conclusion follows from this.

**Conclusion**. To construct a trigonometric spline of one variable, it is enough to construct a trigonometric interpolation polynomial (1) of this variable, and then replace the functions $\cos kx$ and $\sin kx$ spline cos- and sin-functions determined by expressions (4)-(6).

It is clear that this conclusion also applies to the construction of trigonometric splines of several variables.

Let's illustrate the above on the example of constructing a trigonometric spline of two variables.

1. On the XOY plane, consider the region $\Delta = \{0 \le x \le 2\pi, 0 \le y \le 2\pi\}$. As before, let's set a uniform grid $\Delta_{N_1}^{(I_x)} = \left\{ x_j^{(I_x)} \right\}_{j=1}^{N_1}$ by variable $x$, where $I_x$ - grid indicator, $(I_x = 0,1)$, where $x_j^{(0)} = \dfrac{2\pi}{N_1}(j-1)$, $x_j^{(1)} = \dfrac{\pi}{N_1}(2j-1)$, $N_1 = 2n_1 + 1$, $n_1 = 1, 2, \ldots$. Thinking similarly, let's set a uniform grid $\Delta_{N_2}^{(I_y)} = \left\{ y_j^{(I_y)} \right\}_{k=1}^{N_2}$ and byvariable $y$, where $I_y$ - grid indicator, $(I_y = 0,1)$, where $y_k^{(0)} = \dfrac{2\pi}{N_2}(k-1)$, $y_k^{(1)} = \dfrac{\pi}{N_2}(2k-1)$, $N_2 = 2n_2 + 1$, $n_1 = 1, 2, \ldots$.

Since grid indicators must satisfy the condition $I_x = I_y$, in the future we will put $I = I_x = I_y$.

Let it also be in the region $\Delta = \{0 \le x \le 2\pi, 0 \le y \le 2\pi\}$ a function is given $f(x, y)$ two variables. Let's mark $f_{jk}^{(I)} = f\left( x_j^{(I)}, y_k^{(I)} \right)$, $(j = 1, 2, ..., N_1; \; k = 1, 2, ..., N_2)$.

**1.** As is known (see, for example, [12]), the trigonometric interpolation polynomial of two variables has the form

$$T_{n_1 n_2}^{(I)}(x, y) = \sum_{k=0}^{\frac{N_1-1}{2}} \sum_{l=0}^{\frac{N_2-1}{2}} \left( a_{kl}^{(I)} \cos(kx)\cos(ly) + b_{kl}^{(I)} \cos(kx)\sin(ly) + \right. \quad (7)$$
$$\left. + c_{kl}^{(I)} \sin(kx)\cos(ly) + d_{kl}^{(I)} \sin(kx)\sin(ly) \right),$$

where the coefficients $a_{kl}^{(I)}, b_{kl}^{(I)}, c_{kl}^{(I)}, d_{kl}^{(I)}$ are determined by formulas

$$a_{k,l}^{(I)} = K1 \sum_{i=1}^{N_1} \sum_{j=1}^{N_2} f_{i,j}^{(I)} \cos(kx_i^{(I)})\cos(ly_i^{(I)}),$$

$$b_{k,l}^{(I)} = K2 \sum_{i=1}^{N_1} \sum_{j=1}^{N_2} f_{i,j}^{(I)} \cos(kx_i^{(I)})\sin(ly_i^{(I)}), (8)$$

$$c_{k,l}^{(I)} = K3 \sum_{i=1}^{N_1} \sum_{j=1}^{N_2} f_{i,j}^{(I)} \sin(kx_i^{(I)}) \cos(ly_i^{(I)}),$$

$$d_{k,l}^{(I)} = K4 \sum_{i=1}^{N_1} \sum_{j=1}^{N_2} f_{i,j}^{(I)} \sin(kx_i^{(I)}) \sin(ly_i^{(I)}),$$

and the coefficients $K1, K2, K3, K4$ have the form:

$$K1 = \frac{1}{N1}\frac{1}{N2}, K2 = K3 = K4 = 0, \text{ якщо } k = l = 0;$$

$$K1 = \frac{1}{N1}\frac{2}{N2}, K2 = \frac{1}{N1}\frac{2}{N2}, \; K3 = K4 = 0, \text{ якщо } k = 0, l \neq 0;$$

$$K1 = \frac{2}{N1}\frac{1}{N2}, K3 = \frac{2}{N1}\frac{1}{N2}, \; K2 = K4 = 0, \text{ якщо } k \neq 0, l = 0;$$

$$K1 = K2 = K3 = K4 = \frac{2}{N1}\frac{2}{N2}, \text{ якщо } k \neq 0, l \neq 0.$$

Given the conclusion given above, it is easy to construct a trigonometric interpolation spline of two variables, which has the form

$$St_{n_1 n_2}^{(I)}(I_1, I, \nu, r_1, r_2, N_1, N_2, x, y) =$$

$$= \sum_{k=0}^{\frac{N_1-1}{2}} \sum_{l=0}^{\frac{N_2-1}{2}} \left[ a_{kl}^{(I)} Cs_k(I_1, I, \nu, r_1, N_1, x) Cs_l(I_1, I, \nu, r_2, N_2, y) + \right.$$

$$+ b_{kl}^{(I)} Cs_k(I_1, I, \nu, r_1, N_1, x) Ss_l(I_1, I, \nu, r_2, N_2, y) + (9)$$

$$+ c_{kl}^{(I)} Ss_k(I_1, I, \nu, r_1, N_1, x) Cs_l(I_1, I, \nu, r_2, N_2, y) +$$

$$\left. + d_{kl}^{(I)} Ss_k(I_1, I, \nu, r_1, N_1, x) Ss_l(I_1, I, \nu, r_2, N_2, y) \right].$$

Note that the indices in these formulas $I$ interpolation grids are consistent; in addition, it is possible to provide different smoothness properties for each of the arguments.

Reasoning similarly, it is easy to construct a trigonometric interpolation spline of three or more variables; it is clear that the calculation formulas become more complicated.

**2.** It is possible to construct a trigonometric interpolation spline of two variables in another way, using trigonometric fundamental splines. IN[9] it was shown that the system of such splines can be presented in the form

$$ts_k(I_1, I_2, \nu, r, N, x) = \frac{1}{N}\left[ 1 + 2 \sum_{j=1}^{\frac{N-1}{2}} C_j(I_1, I_2, \nu, r, N, k, x) \right], (k = 1, 2, ..., N), (10)$$

where

$$C_j(I_1, I_2, \nu, r, N, k, \mathrm{x}) = \frac{c_j(I_1, I_2, \nu, r, N, k, x)}{h_j(I_1, I_2, \nu, r, N)};$$

$$c_j(I_1, I_2, \nu, r, N, k, x) = \nu(r, j)\cos\left[ j(x - x_k^{(I_2)}) \right] +$$

$$+ \sum_{m=1}^{\infty} (-1)^{mI_1} \left\{ \nu(r, mN - j)\cos\left[ (mN - j)(x - x_k^{(I_2)}) \right] + \nu(r, mN + j)\cos\left[ (mN + j)(x - x_k^{(I_2)}) \right] \right\};$$

$$h_j(I_1, I_2, \nu, r, N) = \nu(r, j) + \sum_{m=1}^{\infty} (-1)^{m(I_1+I_2)}\left[ \nu(r, mN - j) + \nu(r, mN + j) \right].$$

and other designations are the same as before.

Having a system of trigonometric fundamental splines $ts_k(I_1, I_2, \nu, r, N, t)$, $(k = 1, 2, ..., N)$, trigonometric interpolation spline $St(I_1, I_2, \nu, r, N, t)$, which has a stitching grid $\Delta_N^{(I_1)}$ and interpolates the function $f(x)$ in grid nodes $\Delta_N^{(I_2)}$, can be presented in the form

$$St(I_1, I_2, \nu, r, N, t) = \sum_{k=1}^{N} f(x_k^{(I_2)}) ts_k(I_1, I_2, \nu, r, N, t) \quad . \text{(12)}$$

It is clear that the trigonometric interpolation spline of two variables on the domain $\Delta = \{0 \leq x < 2\pi, \ 0 \leq y < 2\pi\}$ can be submitted in the form of

$$St_{n_1 n_2}^{(I)}(I, I, \nu, r_1, r_2, N_1, N_2, x, y) = \sum_{k=1}^{N_1} \sum_{l=1}^{N_2} f_{kl}^{(I)} ts_k(I, I, \nu, r_1, N_1, x) ts_l(I, I, \nu, r_2, N_2, y) \quad . \text{(thirteen)}$$

Note that even in this case, the indices of the interpolation grids are consistent.

Thinking similarly, it is easy to construct trigonometric splines of a larger number of variables.

It is appropriate to emphasize that both methods of constructing a trigonometric interpolation spline of two variables lead to the same results. However, in practice, when constructing trigonometric splines of several variables, it is more convenient to use the second method.

Let's illustrate the above with an example; at the same time, we will limit ourselves to the first option of constructing a trigonometric spline of two variables.

**Example**. Consider the region on the XOY plane $\Delta = \{0 \leq x \leq 2\pi, 0 \leq y \leq 2\pi\}$. As before, let's set a uniform grid $\Delta_{N_1}^{(I_x)} = \{x_j^{(I_x)}\}_{j=1}^{N_1}$ by variable $x$, where $N_1 = 7$, and the grid $\Delta_{N_2}^{(I_y)} = \{y_j^{(I_y)}\}_{k=1}^{N_2}$ on variable $y$ , where $N_2 = 9$. We will also put $I = I_x = I_y$. Let also be a function $f$ takes the value 1 at node 4.5 of the grid $\Delta_{N_1}^{(I)} \times \Delta_{N_2}^{(I)}$ ( $f_{4,5} = 1$ ), and in other nodes of this grid is a function $f$ takes zero values.

Here are some graphs trigonometric interpolation spline of two variables $St_{n_1 n_2}^{(I)}(I, I, \nu, r_1, r_2, N_1, N_2, x, y)$ for some values $r_1 = r_2$ and parameter values $I = 0, 1$.

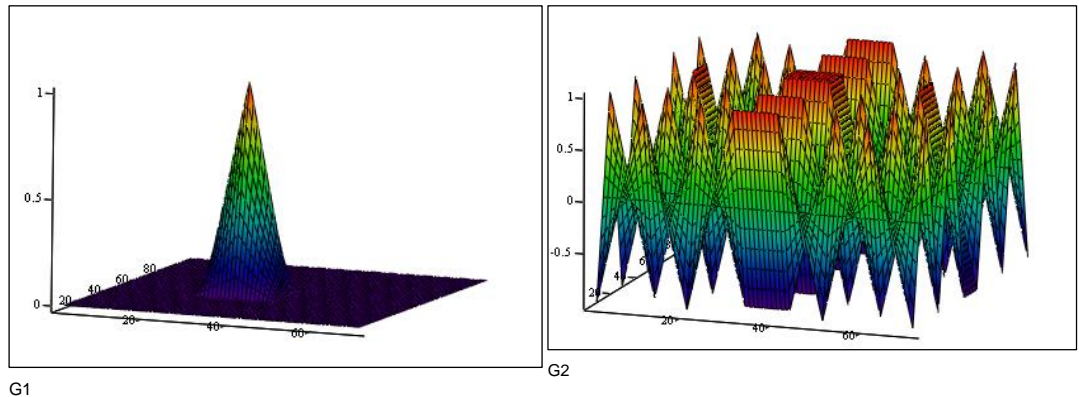

G1                                     G2

Fig. 1. Graphs of the trigonometric interpolation spline of two variables $St_{n_1 n_2}^{(I)}(I, I, \nu, r_1, r_2, N_1, N_2, x, y)$ for values $r_1 = r_2 = 1$ and parameter values $I = 0, 1$.

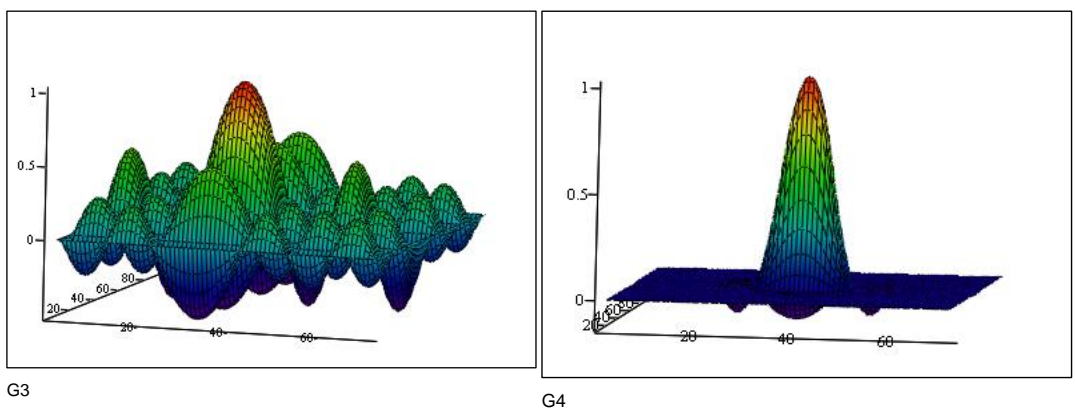

G3    G4

Fig. 1. Graphs of the trigonometric interpolation spline of two variables

$St_{n_1 n_2}^{(I)}(I, I, v, r_1, r_2, N_1, N_2, x, y)$ for values $r_1 = r_2 = 2$

and parameter values $I = 0, 1$.

It follows from the above graphs that with odd parameter values $r_1$ and $r_2$ it is advisable to use splines with parameters $I = 0$; with even parameter values $r_1$ and $r_2$ it is advisable to use splines with parameters $I = 1$. This conclusion can be reached based on other considerations.

## Conclusions

1. A method of constructing trigonometric interpolation splines of two variables of arbitrary degrees on rectangular regions and two types of grids is proposed, in which these splines are obtained based on the interpolation trigonometric polynomial of two variables.
2. A method of constructing trigonometric interpolation splines of two variables of arbitrary degrees on rectangular areas and two types of grids are proposed, in which these splines are obtained as sums of products of fundamental trigonometric splines for each variable.
3. Both methods lead to the same results.
4. Both methods of constructing trigonometric interpolation splines of two variables assume the possibility of constructing such splines with different differential properties for each variable.
5. When constructing trigonometric interpolation splines of several variables of arbitrary degrees on rectangular areas, it is more convenient to use the second construction method, in which these splines are obtained as products of fundamental trigonometric splines for each variable.
6. With odd parameter values $r_1$ and $r_2$ it is advisable to use splines with parameters $I = 0$; with even parameter values $r_1$ and $r_2$ it is advisable to use splines with parameters $I = 1$.
7. The given methods of constructing trigonometric interpolation splines of several variable arbitrary degrees can find wide application in practice.

## List of references